\documentclass[12pt,a4paper]{article}
\usepackage[T1]{fontenc}
\usepackage{mathptmx}
\usepackage{courier}
\usepackage[dvips]{graphicx}
\usepackage{rotating}
\usepackage{psfrag}
\usepackage{amsmath,amssymb}        
\usepackage{latexsym}
\usepackage{url}

\title{A Non-Iterative Transformation Method for an Extended Blasius Problem}
\author{Riccardo Fazio \\
Department of Mathematics, Computer Science\\ Physical Sciences and Earth Sciences,\\
University of Messina \\
Viale F. Stagno D'Alcontres, 31 \\
98166 Messina, Italy \\
E-mail: rfazio@unime.it \\
Home-page: http://mat521.unime.it/fazio}
\date{\today}
\begin{document}
\maketitle
\begin{abstract}
In this paper, we define a non-iterative transformation method for an Extended Blasius Problem.
The original non-iterative transformation method, which is based on scaling invariance properties, was defined for the classical Blasius problem by T\"opfer in 1912.
This method allows us to solve numerically a boundary value problem by solving a related initial value problem and then rescaling the obtained numerical solution. 
In recent years, we have seen applications of the non-iterative transformation method to several problems of interest. 
 
The obtained numerical results are improved by both a mesh refinement strategy and Richardson's extrapolation technique.
In this way, we can be confident that the computed six decimal places are correct.
\end{abstract}
\bigskip

\bigskip


\section{Introduction.}
The problem of determining the steady two-dimensional motion of a
fluid past a flat plate placed edge-ways to the stream was
formulated in general terms, according to the boundary layer theory,
by Prandtl \cite{Prandtl:1904:UFK}, and was investigated in detail by
Blasius \cite{Blasius:1908:GFK}. 
The engineering interest was to calculate the shear
at the plate (skin friction),
which leads to the determination of the viscous drag on the plate, see for instance Schlichting 
\cite{Schlichting:2000:BLT}).
In this contest the celebrated Blasius problem is given by
\begin{align}\label{eq:Blasius} 
& {\displaystyle \frac{d^3 f}{d \eta^3}} + \frac{1}{2} \; f
{\displaystyle \frac{d^{2}f}{d\eta^2}} = 0 \nonumber \\[-1ex]
&\\[-1ex]
& f(0) = {\displaystyle \frac{df}{d\eta}}(0) = 0, \qquad
{\displaystyle \frac{df}{d\eta}}(\eta) \rightarrow 1 \quad \mbox{as}
\quad \eta \rightarrow \infty \ . \nonumber 
\end{align}
This is a boundary value problem (BVP) defined on the semi-infinite interval $[0, \infty)$.
According to Weyl \cite{Weyl:1942:DES}, the unique solution of (\ref{eq:Blasius}) has a positive second order derivative, which is monotone decreasing on $ [0, \infty) $ and approaches to zero as $ \eta $ goes to infinity.
The governing differential equation and the two boundary conditions
at the origin in (\ref{eq:Blasius}) are invariant with respect to the scaling group of transformations
\begin{eqnarray}
& \eta^* = \lambda^{-\alpha}, \qquad f^* = \lambda^{\alpha} f
\label{eq:scaling:Blasius}
\end{eqnarray}
where $ \alpha $ is a nonzero constant, classically $ \alpha = 1/3
$; but here we have used also $ \alpha = 1 $ in order to simplify the
analysis. 
The mentioned invariance property has both analytical
and numerical interest. 
From a numerical viewpoint a non-iterative transformation method (ITM)
reducing the solution of (\ref{eq:Blasius}) to the solution of a related initial value problem (IVP) was defined by T\"opfer \cite{Topfer:1912:BAB}. 
Owing to that transformation, a simple existence and uniqueness Theorem
was given by J. Serrin as reported by Meyer \cite[pp. 104-105]{Meyer:1971:IMF}. 
Let us note here that the mentioned invariance property is essential to the error analysis of the truncated boundary solution due to Rubel \cite{Rubel:1955:EET}, see Fazio \cite{Fazio:2002:SFB}.
Blasius problem was used, recently, by Boyd \cite{Boyd:2008:BFC} as an example were some good analysis, before the computer invention, allowed researchers of the past to solve problems, governed by partial differential equations, that might be otherwise impossible to face.

The application of a non-ITM to the Blasius equation with slip boundary condition, arising within the study of gas and liquid flows at the micro-scale regime \cite{Gad-el-Hak:1999:FMM,Martin:2001:BBL}, was considered already in \cite{Fazio:2009:NTM}.
We applied a non-ITM to the Blasius equation with moving wall considered by Ishak et al. \cite{Ishak:2007:BLM} or surface gasification studied by Emmons \cite{Emmons:1956:FCL} and recently by Lu and Law \cite{Lu:2014:ISB} or slip boundary conditions investigated by Gad-el-Hak \cite{Gad-el-Hak:1999:FMM} or Martin and Boyd \cite{Martin:2001:BBL}, see Fazio \cite{Fazio:2016:NIT} for details.
In particular, we find a way to solve non-iteratively the Sakiadis problem \cite{Sakiadis:1961:BLBa,Sakiadis:1961:BLBb}.
As far as the non-ITM is concerned, a recent review dealing with all the cited problems can be be found in \cite{Fazio:2019:NIT}.

Moreover, T\"opfer's method has been extended to classes of problems in boundary layer theory involving a physical parameter.
This kind of extension was considered first by Na \cite{Na:1970:IVM}, see also NA \cite[Chapters 8-9]{Na:1979:CME}.

Finally, an iterative extension of the transformation method has been introduced, for the numerical solution of free BVPs, by Fazio \cite{Fazio:1990:SNA}. 
This iterative extension has been applied to several problems of interest: 
free boundary problems \cite{Fazio:1990:SNA,Fazio:1997:NTE,Fazio:1998:SAN},
a moving boundary hyperbolic problem \cite{Fazio:1992:MBH},
Homann and Hiemenz problems governed by the Falkner-Skan equation in \cite{Fazio:1994:FSE},
one-dimensional parabolic moving boundary problems \cite{Fazio:2001:ITM}, two variants of the Blasius problem \cite{Fazio:2009:NTM}, namely: a boundary layer problem over moving surfaces, studied first by Klemp and Acrivos \cite{Klemp:1972:MBL}, and a boundary layer problem with slip boundary condition, that has found application in the study of gas and liquid flows at the micro-scale regime \cite{Gad-el-Hak:1999:FMM,Martin:2001:BBL}, parabolic problems on unbounded domains \cite{Fazio:2010:MBF} and,  recently, see \cite{Fazio:2015:ITM}, a further variant of the Blasius problem in boundary layer theory: the so-called Sakiadis problem \cite{Sakiadis:1961:BLBa,Sakiadis:1961:BLBb}.
A recent review dealing with, the derivation and application of, ITM can be be found by the interested reader in \cite{Fazio:2019:ITM}.

\section{Extended Blasius problem}
Our extended Blasius problem is given by,
\begin{align}\label{eq:ExBlasius} 
& {\displaystyle \frac{d^3 f}{d \eta^3}} {\displaystyle \frac{d^{2}f}{d\eta^2}}^{(P-1)} + \frac{1}{2} \; f
{\displaystyle \frac{d^{2}f}{d\eta^2}} = 0 \nonumber \\[-1ex]
&\\[-1ex]
& f(0) = {\displaystyle \frac{df}{d\eta}}(0) = 0, \qquad
{\displaystyle \frac{df}{d\eta}}(\eta) \rightarrow 1 \quad \mbox{as}
\quad \eta \rightarrow \infty \ , \nonumber 
\end{align}
where $P$ verifies the conditions $1 \le P < 2$,  see Schowalter \cite{Schowalter:1960:ABL}, Lee and Ames \cite{Lee:1966:SSN}, Lin and Chern \cite{Lin:1979:LBL}, Kim et al. \cite{Kim:1983:MHT}, or Akcay and Y\"ukselen \cite{Akcay:1999:DRN}.
Liao \cite{Liao:2005:CNP} has found analytically that the extended Blasius problem for $P = 2$ admit an infinite number of solutions anf therefore in his opinion can be considered as a challenge problem for numerical techniques.

\subsection{The non-ITM} 
The applicability of a non-ITM to the Blasius problem (\ref{eq:Blasius}) is a consequence of both: the invariance of the governing differential equation and the two boundary conditions at $\eta = 0$, and the non invariance of the asymptotic boundary condition under the scaling transformation (\ref{eq:scaling:Blasius}).
In order to apply a non-ITM to the BVP (\ref{eq:ExBlasius}) we investigate its invariance with respect to the scaling group
\begin{equation}\label{eq:scaling}
f^* = \lambda f \ , \qquad \eta^* = \lambda^{\delta} \eta \ .   
\end{equation}
We find that the extended Blasius problem (\ref{eq:ExBlasius}) is invariant under (\ref{eq:scaling}) iff
\begin{equation}\label{eq:scaling:condition}
\delta = \frac{2-P}{1-2 P} \ .
\end{equation}
Now, we can integrate the extended Blasius equation in (\ref{eq:ExBlasius}) written in the star variables on $[0, \eta^*_\infty]$, where $\eta^*_\infty$ is a suitable truncated boundary, with initial conditions
\begin{equation}\label{eq:ICs2}
f^*(0) = \frac{df^*}{d\eta^*}(0) = 0 \ , \quad \frac{d^2f^*}{d\eta^{*2}}(0) = 1 \ ,
\end{equation}
in order to compute an approximation $\frac{df^*}{d\eta^*}(\eta^*_\infty)$ for $\frac{df^*}{d\eta^*}(\infty)$ and the corresponding value of $\lambda$ according to the equation
\begin{equation}\label{eq:lambda}
\lambda = \left[ \frac{d f^*}{d \eta^{*}}(\eta^*_\infty) \right]^{1/(1-\delta)} \ .   
\end{equation} 
Once the value of $\lambda$ has been computed by equation (\ref{eq:lambda}), we can find the missed initial condition by the equation
\begin{equation}\label{eq:MIC}
\frac{d^2f}{d\eta^{2}}(0) =  \lambda^{2\delta-1}\frac{d^2f^*}{d\eta^{*2}}(0) \ .
\end{equation}
Moreover, the numerical solution of the original BVP (\ref{eq:ExBlasius}) can be computed by rescaling the solution of the IVP.
In this way we get the solution of a given BVP by solving a related IVP.

\section{Numerical results}
\begin{figure}[!hbt]
	\centering
\psfrag{e}[1][]{$\eta, \eta^*$}  
\psfrag{c}[1][]{$\frac{df^*}{d\eta^*}$,$\frac{d^2f^*}{d{\eta^*}^2}$,$\frac{df}{d\eta}$,$\frac{d^2f}{d\eta^2}$}  
\includegraphics[width=14cm,height=14cm]{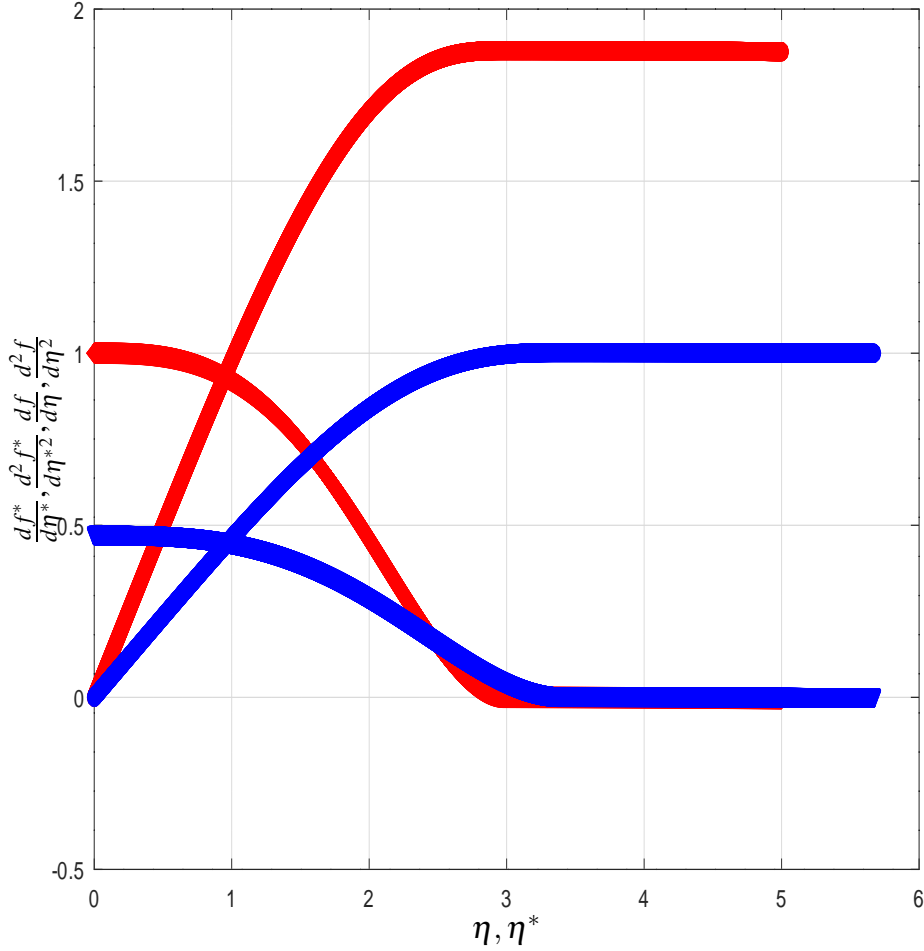}
\caption{Numerical results of the non-ITM for (\ref{eq:ExBlasius}) with $P =3/2$. The starred variables problem and the original problem solution components found after rescaling.}
	\label{fig:ExBlasius}
\end{figure}

In this section we report the numerical results computed with our non-ITM. 
To compute the numerical solution, we used the classical fourth order Runge-Kutta method with constant step size.
For the results shown in this figure we used $\Delta \eta = 0.001$ and $\eta^*_{\infty} = 5$.
Figure \ref{fig:ExBlasius} shows the solution of the extended Blasius problem, describing the behaviour of a boundary layer flow due to a moving flat surface immersed in an otherwise quiescent fluid, corresponding to $P=3/2$.
Let us notice here that, by rescaling, we get $\eta^*_{\infty} < \eta_{\infty}$. 

In table \ref{tab:missingIC} we report the missing initial condition $\frac{d^2f}{d\eta^2}(0)$ obtained by a mesh refinement.
\begin{table}[!hbt]
\caption{Missing initial condition by mesh refinement.}
\vspace{.5cm}
\renewcommand\arraystretch{1.3}
	\centering
		\begin{tabular}{lr@{.}l}
\hline 
{${\Delta \eta}$} &  \multicolumn{2}{c}%
{$\frac{d^2f}{d\eta^2}(0)$} \\[1.2ex]
\hline
0.001         & 0 & 47001559210 \\  
0.0005        & 0 & 46953524356 \\  
0.00025       & 0 & 46929513793 \\  
0.000125      & 0 & 46917509666 \\  
0.0000625     & 0 & 46911507775 \\
0.00003125    & 0 & 46908506833 \\
0.000015625   & 0 & 46907006359 \\
0.0000078125  & 0 & 46906256120 \\
\hline			
		\end{tabular}
	\label{tab:missingIC}
\end{table}
As it is easily seen the decimal digits, starting to the left of the decimal point, tend to be confirmed, in particular, here we have confirmed the first four decimal values.

In order to improve the obtained missing initial condition value, we can apply a Richardson's extrapolation technique, see \cite{Richardson:1910:DAL,Richardson:1927:DAL}.
Let us consider that we would like to compute the value $U = \frac{d^2f}{d\eta^2}(0)$ as $\Delta \eta$ goes to zero, then we can apply the extrapolation formula
\begin{equation}\label{eq:Rextra}
U_{g+1,k+1} = U_{g+1,k} + \frac{U_{g+1,k}-U_{g,k}}{4-1} \ ,
\end{equation}
where $g \in \{0, 1, 2 , \dots , G-1\}$, $k \in \{0, 1, 2, \dots , G-1\}$, here $g$ indicates the $g$-esima computational grid, $k$ the level of extrapolation, and $4$ is the true order of the discretization error for the classical Runge-Kutta method.
The formula (\ref{eq:Rextra}) is asymptotically exact in the limit as $\Delta \eta$ goes to zero if we use uniform grids.
We notice that to obtain each value of $U_{g+1,k+1}$ requires two computed solutions $U$ in two adjacent grids, namely $g+1$ and $g$ at the extrapolation level $k$.
For any $g$, the level $k=0$ represents the numerical solution of $U$ without any extrapolation.
The numerical results obtained applying formula (\ref{eq:Rextra}) to the values reported in table \ref{tab:missingIC} are reported in table \ref{tab:missingICEx}.

\begin{sidewaystable}
\caption{Missing initial condition by Richardson's extrapolation.}
\vspace{.5cm}
\renewcommand\arraystretch{1.3}
	\centering
		\begin{tabular}{lllllllll}
\hline 
{$U_{g,0}$} & {$U_{g,1}$} & {$U_{g,2}$} & {$U_{g,3}$} & {$U_{g,4}$} & {$U_{g,5}$} & {$U_{g,6}$} & {$U_{g,7}$} \\[1.2ex]
\hline
0.47001559210 & & & & & & \\  
0.46953524356 & 0.46937512738 & & & & & & \\  
0.46929513793 & 0.46921510272 & 0.46916176117 & & & & & \\  
0.46917509666 & 0.46913508290 & 0.46910840963 & 0.46909062578 & & & & \\  
0.46911507775 & 0.46909507145 & 0.46908173430 & 0.46907284252 & 0.46906691477 & & & \\  
0.46908506833 & 0.46907506519 & 0.46906839644 & 0.46906395049 & 0.46906098648 & 0.46905901038 & & \\
0.46907006359 & 0.46906506201 & 0.46906172762 & 0.46905950468 & 0.46905802274 & 0.46905703482 & 0.46905637630 &  \\
0.46906256120 & 0.46906006040 & 0.46905839320 & 0.46905728173 & 0.46905654074 & 0.46905604674 & 0.46905571739 & 0.46905549774  \\
\hline			
		\end{tabular}
	\label{tab:missingICEx}
\end{sidewaystable}

Naturally, we can stop the mesh refinement as soon as $U_{g+1,k} = U_{g,k}$ or Richardon's extrapolation when $U_{g,k+1} = U_{g,k}$.
It is evident, from the data reported in table \ref{tab:missingICEx}, that we have achieved six decimal places of agreement.

As mentioned before, the case $P=1$ is the Blasius problem (\ref{eq:Blasius}).
In this case our non-ITM reduces to the original method defined by T\"opfer \cite{Topfer:1912:BAB}.
In that case, the computed skin friction coefficient value, namely $0.332057336215$, obtained with $\Delta \eta = 0.001$ and $\eta^*_{\infty} = 10$, is in good agreement with the values available in literature, see for instance the value $0.332057336215$ computed by Fazio \cite{Fazio:1992:BPF} by a free boundary formulation of the Blasius problem or the value $0.33205733621519630$ computed by Boyd \cite{Boyd:1999:BFC} who believes all the decimal digits to be correct.

\section{Concluding remarks.}
The main contribution of this paper is the extension of the non-ITM, proposed by T\"opfer \cite{Topfer:1912:BAB} and defined for the numerical solution of the celebrated Blasius problem \cite{Blasius:1908:GFK}, to an extended Blasius problem.
This method allows us to solve numerically the extended Blasius problem by solving a related initial value problem and then rescaling the obtained numerical solution. 
The obtained numerical results, have been improved both by a mesh refinement and using the Richardson's extrapolation technique.
In this way, we can be confident that the computed six decimal place are correct.

\vspace{1.5cm}

\noindent {\bf Acknowledgement.} {The research of this work was 
partially supported by the University of Messina and by the GNCS of INDAM.}


\end{document}